\documentclass{article}
\pdfoutput=1 

\usepackage{arxiv}

\DeclareMathAlphabet{\pazocal}{OMS}{zplm}{m}{n}
\usepackage{mathrsfs}
\usepackage{calligra}
\usepackage{amssymb,amsmath,amsfonts,amsthm}
\usepackage{bm}
\usepackage{dsfont}
\usepackage{pifont}
\usepackage{mathtools}

\usepackage[utf8]{inputenc}
\usepackage[T1]{fontenc} 	
\usepackage{hyperref}
\usepackage{url}            
\usepackage{booktabs}       
\usepackage{titlesec}
\usepackage{graphicx}
\usepackage{enumerate}
\usepackage{enumitem}
\usepackage[super]{nth}
\usepackage{xcolor}
\usepackage{tabularx}
\usepackage{authblk}
\usepackage{ragged2e}
\definecolor{darkpastelgreen}{rgb}{0.01,0.75,0.24}

\theoremstyle{plain}
\newtheorem{theo}{Theorem}[section] 
\newtheorem{lemma}[theo]{Lemma} 

\newtheorem{coroll}[theo]{Corollary}

\theoremstyle{definition}
 
\newtheorem{exmp}[theo]{Example} 
\newtheorem{com}[theo]{Remark}

\titleformat{\section}[hang]
  {\bfseries}
  {\thesection.}
  {1ex}
  {}
\titlespacing{\section}{1.5pt}{0.2cm}{0.2cm}

\titleformat{\subsection}[hang]
  {\bfseries}
  {\thesubsection.}
  {1ex}
  {}
\titlespacing{\subsection}{1.5pt}{0.2cm}{0.2cm}

\newcommand{\1}{\partial}
\newcommand{\w}{\omega}

\newcommand{\grad}{\text{grad}_{\textsl{g}}\hspace{0.05cm}}

\newcommand{\lna}{\text{log}\hspace{0.05cm} }

\newcommand{\pzero}{\textit{p}\hspace{-0.02cm}_{\textit{0}}}
\newcommand{\expo}{\text{exp}}
\newcommand{\hit}{\textit{h}}

\newcommand{\emef}{\textit{m}_f}
\newcommand{\eme}{\textit{m}}
\newcommand{\emeh}{\textit{m}_{\hspace{-0.03cm}\textbf{\tiny{0}}}}

\newcommand{\has}{\mathrm{H}}
\newcommand{\hasf}{\mathrm{H}_f}

\newcommand{\Ricc}{\mathrm{R}\hspace{-0.03cm}\mathrm{i}\hspace{-0.045cm}\mathrm{c}}
\newcommand{\Riccfn}{\mathrm{R}\hspace{-0.03cm}\mathrm{i}\hspace{-0.045cm}\mathrm{c}_{f}^{N}}
\newcommand{\Riccf}{\mathrm{R}\hspace{-0.03cm}\mathrm{i}\hspace{-0.045cm}\mathrm{c}_{f}}

\renewcommand{\shorttitle}{\textit{Remark on a geometric inequality in weighted manifolds}}

\hypersetup{
pdftitle={Somemathtitle},
pdfsubject={Comparison Geometry, Differential Geometry},
pdfauthor={Adam Rudnik},
pdfkeywords={Asympotically nonnegative curvature, Asymptotic volume ratio, Willmore-type inequality},
}

\begin{document}
\makeatletter
\newcommand{\address}[1]{\gdef\@address{#1}}
\newcommand{\email}[1]{\gdef\@email{\url{#1}}}
\newcommand{\@endstuff}{\par\vspace{\baselineskip}\noindent\small
\begin{tabular}{@{}l}\scshape\@address\\\textit{E-mail address:} \@email 
\end{tabular}}
\AtEndDocument{\@endstuff}
\makeatother

\title{%
	\boldmath REMARK ON A GEOMETRIC INEQUALITY FOR CLOSED HYPERSURFACES IN WEIGHTED MANIFOLDS	\\ \vspace{0.4cm}
	\large\textit{To the memory of Celso Viana}
	}

\author{\href{https://orcid.org/0000-0003-3181-1466}{\includegraphics[scale=0.06]{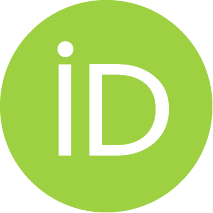}\hspace{1mm}ADAM RUDNIK}}

\address{\parbox{\linewidth}{Department of Mathematics, Institute of Exact Sciences (ICEx), Federal University of Minas Gerais,\\ Belo Horizonte, MG 31270901, Brazil}}
\email{adamrudk@ufmg.br}

\date{}
\maketitle

\vspace{-0.5cm}
\begin{abstract}\footnotesize
In this paper we consider noncompact smooth metric measure spaces $(\textsl{M}, \textsl{g},e^{-f}\textit{dvol}_{\textsl{g}})$ of nonnegative Bakry-Émery Ricci curvature, i.e. $\Ricc + D^{2}f - \frac{1}{N}df \otimes df \geq 0$, for $0< N \leq \infty$, in order to obtain geometric inequalities for the boundary of a given open and bounded set $\Omega\subset \textsl{M}$, with regular boundary $\1 \Omega$. Our inequalities are sharp for both the cases $N< \infty$ and $N= \infty$, provided that the underlying ambient space has large weighted volume growth. The rigidity obtained for the $N=\infty$ case holds true precisely when $\textsl{M}\setminus \Omega$ is isometric to a twisted product metric and, as such, is a generalization of the Willmore-type inequality for nonnegative Ricci curvature from \cite{agostiniani2020sharp} to the context of weighted manifolds.
\end{abstract}


\section{Introduction}
In this work, we consider noncompact $n$-dimensional smooth metric measure spaces $(\textsl{M}, \textsl{g},e^{-f}\textit{dvol}_{\textsl{g}})$ with nonnegative Bakry-Émery Ricci curvature, that is

\begin{flalign}
\label{general_hyp}
\hspace{0.8cm}
\Ricc + D^{2}f - \frac{1}{N}df \otimes df & \geq 0, \hspace{0.2cm} \text{where} \hspace{0.2cm} 0< N \leq \infty , &
\end{flalign}
$f$ is a smooth real function on \textsl{M}, and $\Ricc$ is the Ricci curvature of the pair $(\textsl{M}, \textsl{g})$. Formally, a smooth metric measure space (abbreviated SMMS) is a triple $(\textsl{M}, \textsl{g},e^{-f}\textit{dvol}_{\textsl{g}} )$, where $(\textsl{M}, \textsl{g})$ is a complete Riemannian manifold together with a weighted volume form $e^{-f}\textit{dvol}_{\textsl{g}}$, where $\textit{dvol}_{\textsl{g}}$ is the standard volume form on $\textsl{M}$ induced by $\textsl{g}$. Given a SMMS $(\textsl{M}, \textsl{g},e^{-f}\textit{dvol}_{\textsl{g}} )$ satisfying the lower bound (\ref{general_hyp}), and $\pzero\in \textsl{M}$, relative volume comparison theorems (see \cite{wei2009comparison}, \cite{qian1997estimates}) say that the functions 

\begin{flalign}
\label{rel_vol_grow}
\hspace{0.8cm}
r\mapsto & \frac{\textit{vol}_{f}(B_{\pzero}(r))}{\w _{n+N} r^{n+N}}, \hspace{0.2cm}\text{if}\hspace{0.2cm}  N <\infty, \hspace{0.2cm} \text{and} \hspace{0.2cm}
 r\mapsto \frac{\textit{vol}_{f}(B_{\pzero}(r))}{\w _{n} r^{n}}, \hspace{0.2cm}\text{if}\hspace{0.2cm}  N = \infty , &
\end{flalign}
are nonincreasing, where $\w _{k}$ is the volume of a $k$-dimensional Euclidean unit ball. In the case $N=\infty$, we further assume $\langle \nabla f, \nabla \rho \rangle\geq 0$, where $\rho = \text{dist}_{\textsl{g}}(\cdot,\pzero)$ is the distance to $\pzero$. Then, we define the \textit{$f\textit{-}$Asymptotic Volume Ratio} of $\textsl{g}$, denoted here by $f\textit{-}AVR(\textsl{g})$, as being the limit of the quotient in (\ref{rel_vol_grow}) as $r\rightarrow \infty$. The SMMS $(\textsl{M}, \textsl{g},e^{-f}\textit{dvol}_{\textsl{g}} )$ is said to have  \textit{large weighted volume growth} if $f\textit{-}AVR(\textsl{g})>0$, for some $\pzero \in \textsl{M}$. When $f$ is constant, (\ref{general_hyp}) reduces to $\Ricc \geq 0$ while the $f\textit{-}$asymptotic volume ratio reduces to the well known classical asymptotic volume ratio of the pair $(\textsl{M}, \textsl{g})$ and large weighted volume growth reduces to large (or Euclidean) volume growth, in the classical sense. Metrics with nonnegative Ricci curvature were extensively studied in the past several decades (\cite{shen1993volume}, \cite{sha1997complete}, \cite{sha1989examples}, \cite{perelman1994proof}) and special attention has been given to spaces with large volume growth (\cite{cheeger1971splitting}, \cite{li1986large}, \cite{anderson1990topology}, \cite{shen1996complete}, \cite{carmo2012ricci}). Therefore, condition (\ref{general_hyp}) is a rather natural generalization of nonnegative Ricci curvature.\vspace{0.2cm}

Weighted volume measures arise naturally from the study of conformal deformations of a Riemannian metric (\cite{chang2006conformal}), also as smooth collapsed limits of spaces with lower bounds on Ricci curvature under measured Hausdorff-Gromov convergence (\cite{fukaya1987collapsing}) and in relation to optimal transport theory (\cite{lott2009ricci}). For instance, given two conformally related Riemannian metrics, it is desirable to establish geometric and topological results relatively to the new metric using data associated to the original metric. To accomplish this, it is interesting to have a concept of curvature associated to a weighted Laplacian. Such a concept has been introduced by Bakry and Émery \cite{bakry2006diffusions}. The Bakry-Émery Ricci tensor gives an analogue of the Ricci tensor for a Riemannian manifold with smooth measure. Therefore, it is interesting to investigate what geometric and topological results for the Ricci tensor extend to the Bakry-Émery Ricci tensor. The Bakry-Émery Ricci tensor has an extension for diffusion operators (\cite{bakry2006diffusions}) and nonsmooth metric measure spaces. Moreover, the study of lower bounds on the Bakry-Émery Ricci tensor leads to geometrical and topological information on gradient Ricci solitons (\cite{munteanu2011smooth}, \cite{munteanu2012analysis}). Furthermore, for $N$ positive integer, the equation $\Riccfn = \lambda \textsl{g}$ corresponds to a warped product Einstein metric on $\textsl{M}\times_{e^{-\frac{f}{N}}}\textsl{F}^{N}$, see \cite{case2011rigidity},  \cite{he2012classification}.

In order to systematize our study and easily distinguish between the $N<\infty$ and the $N=\infty$ case, we introduce the following notation. The \textit{N-Bakry-Émery Ricci tensor} is defined by
\begin{flalign}
\label{finite_bakry-emery_Ricci-T}
\hspace{0.8cm}
\Riccfn & \doteq \Ricc + D^{2}f - \frac{1}{N}df \otimes df, &
\end{flalign}
for positive $N >0$, and the \textit{$\infty$-Bakry-Émery Ricci tensor} is defined by 

\begin{flalign}
\label{infinite_bakry-emery_Ricci-T}
\hspace{0.8cm}
\Riccf & \doteq \Ricc + D^{2}f. &
\end{flalign}
In the past few years, geometric inequalities of Willmore-type have received a large amount of attention \cite{agostiniani2020monotonicity}, \cite{borghini2023comparison}, \cite{celso23iso}. For spaces with nonnegative Ricci curvature, a breakthrough in the subject is certainly the work of Agostiniani, Fagagnolo and Mazzieri 	\cite{agostiniani2020sharp}, where the authors obtained a Willmore-type inequality that is sharp for spaces with large volume growth. Nonetheless, it was the simplification obtained by Wang \cite{wang2022} of the very same inequality that paved the way for other works on the subject, e.g. \cite{rudnik2023sharp}. In this context, it is natural to consider the analogous rigidity problem under nonnegative Bakry-Émery Ricci curvature, that is, the existence of such geometric inequalities under the general curvature condition (\ref{general_hyp}). In \cite{wu2025willmore}, Wu and Wu pursue exactly that investigation. However, for a Willmore-type inequality under nonnegative Bakry-Émery Ricci curvature, the question if one can obtain the same qualitative sharpness, as in the case for nonnegative Ricci curvature, remains open. The main goal of this paper is to give an affirmative partial answer to this issue. \vspace{0.2cm}

The purpose of the present work is to obtain two geometric inequalities for a class of closed hypersurfaces $\Sigma$ in noncompact manifolds $\textsl{M}$ satisfying the lower bound (\ref{general_hyp}), for some function $f$ on $\textsl{M}$, which are sharp for spaces with large weighted volume growth. The rigidity achieved here, i.e., the equality case of these geometric inequalities, improve the previous results of Wu and Wu, where the authors obtained similar conclusions under the stronger hypothesis that the weighted mean curvature of the hypersurface $\Sigma$ is constant. More precisely, Theorem \ref{main_theo_N_bakry-emery} ahead extends Theorem 1.4 from \cite{wu2025willmore}. Moreover, the rigidity in Theorem \ref{main_theo_N_bakry-emery} is similar to the rigidity in the main theorem in \cite{agostiniani2020sharp} in: \textit{(a)} both yield that the boundary $\1\Omega$ is an umbilic hypersurface of constant mean curvature; and \textit{(b)} the truncated cone over $\1\Omega$ is the only geometry which gives equality. The $N =\infty$ case is handled analogously. We show that the exact model geometry, corresponding to equality, is actually a twisted product metric over the referred hypersurface, extending and generalizing the Willmore-type inequality from \cite{agostiniani2020sharp} to manifolds with density  (see Theorem \ref{main_theo_infty_bakry-emery} for details). In what follows, we assume $N\geq 1$ to be an integer. Our first main theorem is the following. \vspace{0.2cm}


\begin{theo} $\big[$\textit{Willmore-type inequality for nonnegative N-Bakry-Émery} $\big]$ \newline
\label{main_theo_N_bakry-emery}
Let $(\textsl{M},\textsl{g}, e^{-f}\textit{dvol}_{\textsl{g}})$ be a noncompact, $n$-dimensional, SMMS with nonnegative N-Bakry-Émery Ricci curvature, $\Riccfn \geq 0$. Let $\Omega$ be an open and bounded set with smooth boundary $\Sigma = \partial \Omega$, whose mean curvature is $\has$. Then

\begin{flalign}
\label{willmorex_type_N_bakry-emery}
\hspace{0.8cm}
\int_{\Sigma}\bigg \vert \frac{ \hasf(x)}{n+N-1} \bigg\vert^{n+N-1} e^{-f(x)}d\sigma (x) & \geq \big\vert\mathbb{S}^{n+N-1} \big\vert f \textit{-}AVR(\textsl{g}), &
\end{flalign}
where $\hasf = \has - \langle \nu, \grad f \rangle$ is the weighted mean curvature of $\Sigma$ and $\nu$ is the outward unit normal. Furthermore, if $f \textit{-}AVR(\textsl{g})>0$ and $\Sigma$ is connected then equality in (\ref{willmorex_type_N_bakry-emery}) holds iff $\textsl{M} \setminus \Omega$ is isometric to

\begin{flalign}
\label{warped_product_split_N_bakry-emery}
\hspace{0.8cm}
\bigg([r_0, \infty )\times \Sigma, & dr^2 + (r/r_0)^2 g_{\Sigma}\bigg), \hspace{0.3cm} \text{where} \hspace{0.3cm} r_0 = \frac{n-1}{\has} &
\end{flalign}
and outside $\Omega$, f is given by $f(\exp_x r\nu_x) = \lna\big( 1+ \frac{\has r}{n-1}\big)^{-N} + f(x)$, for all $x \in \Sigma$. In particular, $\Sigma$ is a totally umbilic hypersurface with constant (weighted) mean curvature. 
\vspace{0.2cm}
\end{theo}

\begin{com}
We point out that the quantity $f\textit{-}AVR(\textsl{g})$ seems to not depend on the point $\pzero$ in \textsl{M}, for $N< \infty$ (see \cite{Johne2021Sobolev}). Moreover, we expect that Theorem \ref{main_theo_N_bakry-emery} extends from integers $N\geq 1$ to real numbers $N>0$. To obtain this, one should use the volume formula $V_x = \frac{\pi^{x/2}}{\Gamma(x/2 +1)}$ for the unit $x$-ball, in terms of the Gamma function $\Gamma$, so that $A_{x-1} = \frac{2\pi^{x/2}}{\Gamma(x/2)}$ is the area of the unit $(x-1)$-sphere.
\end{com}

\begin{com}
\label{topol_rem_N<infty_will}
Theorem \ref{main_theo_N_bakry-emery} is in fact a byproduct of a more general result which we now explain. Letting the connectedness assumption on $\1\Omega$ aside, we obtain the following. If $(\textsl{M},\textsl{g}, e^{-f}\textit{dvol}_{\textsl{g}})$ has large weighted volume growth and equality in (\ref{willmorex_type_N_bakry-emery}) holds then we consider the set $\Sigma^+$, defined as the subset of $\Sigma$ where the weighted mean curvature is positive. We prove that $\Sigma^+$ is a totally umbilical submanifold with constant (weighted) mean curvature on its connected components and obtain the preceding warped product decomposition (\ref{warped_product_split_N_bakry-emery}) on the portion of $\textsl{M}\setminus \Omega$ equal to the image of $\textsl{N}\Sigma^+$ under the normal exponential map. See the proof of Theorem \ref{main_theo_N_bakry-emery} for details.
\end{com}
Our corollary may be regarded as a Kasue nonexistence type result (cf. \cite{kasue1983}) for spaces with nonnegative $N$-Bakry-Émery Ricci tensor and large weighted volume growth. It is an immediate consequence from Theorem \ref{main_theo_N_bakry-emery}. \vspace{0.2cm}

\begin{coroll} $\big[$\textit{Non-existence of closed, $f$-minimal hypersurfaces I} $\big]$ \newline
\label{non_existence_N_bakry-emery}
If $(\textsl{M},\textsl{g}, e^{-f}\textit{dvol}_{\textsl{g}})$ is a noncompact, $n$-dimensional, SMMS with $\Riccfn \geq 0$ then there is no $f$-minimal hypersurface $\Sigma$ realized as the boundary of an open and bounded set, provided $f \textit{-}AVR(\textsl{g})>0$. 
\end{coroll}

Before we continue to state our next theorem, let us recall the notion of twisted products (cf. \cite{chen2017differential}). Let $(\textsl{B},\textsl{g}_{\textsl{B}})$ and $(\textsl{F}, \textsl{g}_{\textsl{F}})$ be Riemannian manifolds and let $f: \textsl{B}\times \textsl{F}\to \mathbb{R}$ be a smooth positive function. Then the twisted product of $\textsl{B}$ and $\textsl{F}$ with twisting function $f$ is defined to be the product manifold $\textsl{B} \times \textsl{F}$ provided with the metric tensor $\textsl{g}_{\textsl{B}} + f^2\textsl{g}_{\textsl{F}}$. 

\begin{theo} $\big[$\textit{Willmore-type inequality for nonnegative $\infty$-Bakry-Émery} $\big]$ \newline
\label{main_theo_infty_bakry-emery}
Let $(\textsl{M},\textsl{g}, e^{-f}d \textit{vol}_{\textsl{g}})$ be a noncompact SMMS with nonnegative $\infty$-Bakry-Émery Ricci curvature. Let $\Omega$ be an open and relatively compact set with smooth boundary $\Sigma = \partial \Omega$, whose mean curvature is $\has$. Let $\pzero \in \Omega$ and assume $\frac{\1 f}{\1 r} \geq 0$ and $\frac{\1 f}{\1 \rho} \geq 0$, where $r = \text{dist}_{\textsl{g}}(\cdot , \Sigma)$, is the (signed) distance to $\Sigma$, and $\rho = \text{dist}_{\textsl{g}}(\cdot , \pzero)$. If $\hasf \geq 0$ then

\begin{flalign}
\label{willmorex_type_infty_bakry-emery}
\hspace{0.8cm}
\int_{\Sigma}\bigg ( \frac{ \hasf(x)}{n-1} \bigg)^{n-1} e^{-f(x)}d\sigma (x) & \geq f \textit{-}AVR(\textsl{g}) \vert\mathbb{S}^{n-1} \vert &
\end{flalign}
where $\hasf = \has - \nu (f)$ is the weighted mean curvature of $\Sigma$ and $\nu$ the outward unit normal. Furthermore, if $f \textit{-}AVR(\textsl{g})>0$ and $\hasf >0$ on $\Sigma$ then equality holds iff $\textsl{M}\setminus \Omega$ is isometric to the twisted product metric
\begin{flalign}
\label{twisted_product_m}
\hspace{0.8cm}
\bigg( [0, \infty)\times \Sigma,& dr^2 + \bigg( 1+\frac{\has(x)}{n-1}r \bigg)^2g_{\Sigma}\bigg).&
\end{flalign}
In particular, $\Sigma$ is a totally umbilic hypersurface and $f = \lna H^{n-2}$ outside $\Omega$. 
\end{theo}

\begin{com} If $\Riccf \geq 0$ and $\1_{\rho}f \geq 0$ then the $f\textit{-}$Asymptotic Volume Ratio may, in general, depend on $\pzero$. For this reason, we assume $\pzero \in \Omega$ so that $\textit{vol}_{f}(B_{\pzero}(r)) \leq \textit{vol}_{f}(\mathcal{T}_{\Omega}(r))$, where $\mathcal{T}_{\Omega}(r)$ is the geodesic tube of radius $r$ about $\Omega$ (see Section \ref{proofs_section} for the definition). This allows one to control the relative $f\textit{-}$volume growth of metric balls $B_{\pzero}(r)$, see \ref{prf_ineq_eight} for details. 
\end{com}

\begin{com}
If the mean curvature $\has>0$ in the model (\ref{twisted_product_m}) is constant, then we recover the model geometry in the main theorem from Wu and Wu (\cite[Theorem~1.1,~(a)]{wu2025willmore}) relatively to the curvature condition $\Riccf \geq 0$. \vspace{0.2cm}
\end{com}

\begin{com}
\label{anything}
Likewise Remark \ref{topol_rem_N<infty_will}, in the proof of the rigid case in Theorem \ref{main_theo_infty_bakry-emery}, before proceeding to assume that the weighted mean curvature is positive on the whole of $\Sigma$, we first obtain the twisted product metric decomposition in (\ref{twisted_product_m}) on the subset of $\Sigma$ where the weighted mean curvature is positive and then we globalize the result, assuming $\hasf >0$ on $\Sigma$. \vspace{0.2cm}
\end{com}

\begin{coroll} $\big[$\textit{Non-existence of closed, $f$-minimal hypersurfaces II} $\big]$ \newline
\label{non_existence_infty_bakry-emery}
Under the same conditions of Theorem \ref{main_theo_infty_bakry-emery}, if $f \textit{-}AVR(\textsl{g})>0$ then there is no $f$-minimal hypersurface in $\textsl{M}$ realized as the boundary of an open and relatively compact set. \vspace{0.2cm}
\end{coroll}

\section{Preliminaries}
In this section we establish basic notions, notations and mean curvature comparisons, which are the basic ingredients to obtain the type of volume comparisons we wish. After this, a word on the corresponding asymptotic volume ratio is in order.

\subsection{Notions and notations}
\label{not_not}

Let $\Omega$ be an open, relatively compact set in $\textsl{M}$ with smooth boundary $\1\Omega =\Sigma$. The letter $x$ denotes the variable in $\Sigma$ and
$w$ the variable in $\textsl{M}$. Define the signed distance to the hypersurface $\Sigma$ by 

\begin{flalign}
\label{distance_to_hyper}
\hspace{0.8cm}
r(w)& = \text{dist}_{\textsl{g}}(w, \Sigma), \hspace{0.2cm} w \in \textsl{M}\setminus \Omega. &
\end{flalign}
Denote by $\nu :\Sigma \to T\textsl{M}$ the unit, outward-pointing (relatively to $\Omega$), smooth and normal to $\Sigma$ vector field, and by $N \Sigma$ the \textit{normal bundle} of $\Sigma$. The \textit{normal exponential map} is the restriction of the exponential map from $T\textsl{M}$ to $N \Sigma$. If $SN\Sigma$ denotes those vectors in $N\Sigma$ of norm $\mathds{1}$ then, given $x \in \Sigma$, define the \textit{cut time of $(x, \nu_x)$} by 

\begin{flalign}
\label{cut_time}
\hspace{0.8cm} \tau(\nu_x)& \doteq \text{sup } \{ b>0: \hspace{0.1cm} r(\gamma_{\nu_x}(b)) = b  \}&
\end{flalign}
where the curve $\gamma_{\nu_x}$ is defined as $[0, \infty)\ni t \mapsto \expo_{x}t\nu_x$. We write, unambiguously, $\tau(x) =  \tau(\nu_x)$. Denote the \textit{cut locus of $\Sigma$} by $C( \Sigma)$, so that the normal exponential map from $\{ t\nu_x : x \in \Sigma \hspace{0.2cm} \text{and} \hspace{0.2cm} 0 \leq t < \tau (x)\}$ to $(\textsl{M}\setminus \Omega) \setminus C(\Sigma)$ is a diffeomorphism (see \cite{petersen2006riemannian}, pp. 139 \textit{--} 141). As a result, there exists a diffeomorphism 

\begin{flalign}
\label{diffeo_Phi}
\hspace{0.8cm}
\Phi &: E  \doteq \{ (x,r) \in \Sigma \times [0, \infty) : r < \tau(x)  \}  \to (\textsl{M}\setminus \Omega) \setminus C( \Sigma), \hspace{0.2cm} \Phi(x, r) = \expo_{x}r\nu_x.  &
\end{flalign}
Next, we define the real function $\mathcal{A}$ on $E$ by $(\Phi^{*}d \textit{vol}_{\textsl{g}})_{(x, r)} = \mathcal{A}(x,r) d\sigma \wedge dr$ where $d\sigma$ is the volume element for $\Sigma$ induced by $\textsl{g}$. Therefore, if $f\in C^{\infty}(\textsl{M})$ is a function

\begin{flalign}
\label{coefficient}
\hspace{0.8cm}
\big(\Phi^{*}( e^{-f}d \textit{vol}_{\textsl{g}})\big)_{(x, r)} & = \mathcal{A}_{f}(x,t) d\sigma \wedge dr, \hspace{0.2cm} \text{where} \hspace{0.2cm} \mathcal{A}_{f} \doteq e^{f \circ \Phi}\mathcal{A}. &
\end{flalign}
For $x \in \Sigma$ and any $X,Y\in T_x\Sigma$, the shape operator $A_{\xi}$ of $\Sigma$ is defined by

\begin{flalign*}
\hspace{0.8cm}
\langle A_{\xi} X,Y \rangle & = \langle \nabla_{X}Y, \xi \rangle, \hspace{0.2cm} \xi\in N_x\Sigma. &
\end{flalign*}
The \textit{mean curvature vector $\vec{H}$ of $\Sigma$} at $x$ is $\vec{H}(x) \doteq \sum_{i=1}^{n-1}\langle \nabla_{X_i}X_i, \nu \rangle \nu$, in terms of an orthonormal basis $X_1, ..., X_{n-1}$ of $T_{x}\Sigma$. The mean curvature $\has$ of $\Sigma$ is the trace of $A$ so that $\vec{H} = -\has \nu$ and the weighted mean curvature of $\Sigma$, relatively to $f\in C^\infty(M)$, is $\hasf = \has - \nu(f)$. 
\vspace{0.1cm}

\subsection{Mean curvature and volume comparison}
\label{mean_curv_comp_bakry-emery}
An approach to comparison geometry based on Ricci curvature is given by the \textit{Bochner formula}

\begin{flalign*}
\hspace{0.8cm}
\frac{1}{2}\Delta \vert \nabla u \vert^2 & = \vert D^{2}u \vert^2 + \langle \nabla u, \nabla \Delta u \rangle +\Ricc (\nabla u, \nabla u), \hspace{0.2cm} u \in C^{\infty}(\textsl{M}). &
\end{flalign*}
Relatively to the weighted measure $e^{-f}d \textit{vol}_{\textsl{g}}$, where $f$ is any smooth function on $\textsl{M}$, the intrinsic self-adjoint weighted laplacian is defined by

\begin{flalign*}
\hspace{0.8cm}
\Delta_f u & = \Delta u - \langle \nabla u, \nabla f \rangle, \hspace{0.2cm} u \in C^{\infty}(\textsl{M}). &
\end{flalign*}
It is easy to deduce the Bochner identity for the weighted laplacian $\Delta_f$ in conjunction with the Bakry-Émery Ricci tensors (\ref{finite_bakry-emery_Ricci-T}) and (\ref{infinite_bakry-emery_Ricci-T}) (refer to \cite{wei2009comparison}, \cite{wei2007comparison} for these formulas). For any $u \in C^{\infty}(\textsl{M})$, we have

\begin{flalign}
\label{mean_curv_infty-bakry}
\hspace{0.8cm}
\frac{1}{2}\Delta_f \vert \nabla u \vert^2 & = \vert D^{2}u \vert^2 + \langle \nabla u, \nabla \Delta_f u \rangle + \Riccf (\nabla u, \nabla u) \notag\\
  & \geq \frac{(\Delta u )^2}{n} + \langle \nabla u, \nabla \Delta_f u \rangle + \Riccf (\nabla u, \nabla u), &
\end{flalign}
and the inequality (\ref{mean_curv_infty-bakry}) provides us with the Bochner inequality needed for the $N=\infty$ case. For the $N$ finite case, using the inequality 

\begin{flalign*}
\hspace{0.8cm}
\frac{(\Delta u )^2}{n} + \frac{1}{N}\langle \nabla u, \nabla f \rangle^2 & \geq \frac{(\Delta_f u)^2}{n+N}&
\end{flalign*}
together with (\ref{mean_curv_infty-bakry}), we obtain 

\begin{flalign}
\label{bochner_f_N_bakry}
\hspace{0.8cm}
\frac{1}{2}\Delta_f \vert \nabla u \vert^2 & = \vert D^{2}u \vert^2 + \langle \nabla u, \nabla \Delta_f u \rangle + \Riccf (\nabla u,\nabla u) \notag \\
	& \geq \frac{(\Delta_f u)^2}{n+N} + \langle \nabla u, \nabla \Delta_f u \rangle + \Riccfn (\nabla u, \nabla u).  &
\end{flalign}
We now turn our attention to mean curvature comparison. Recall that the mean curvature measures the relative rate of change of the volume element of $(\text{dist}_{\textsl{g}}(\cdot , \Sigma))^{-1}(t)$. With respect to the measure $e^{-f}d \textit{vol}_{\textsl{g}}$, the weighted mean curvature of the level sets is $\emef = \eme - \frac{\1}{\1 r}(f)$ where $\eme $ is the mean curvature relatively to the inward pointing normal $-\grad r$, so that, $\emef(0) = \hasf$.  \\

Fix any $x\in \Sigma$. The corresponding quantity in the model space we are comparing $\emef$ with is defined by 
\begin{flalign}
\label{mean_curv_model}
\hspace{0.8cm}
\emeh^{k}(r) & \doteq \frac{(k-1) \hasf(x)}{k-1 + \hasf(x) r}, \hspace{0.2cm} \hasf(x) = \has(x) - \nu_x(f), &
\end{flalign}
where $k=n+N$ or $k = n$ depending on whether $N < \infty$ or $N = \infty$, respectively. In general, we will omit the superscript $k$ from $\emeh^{k}$ and it is supposed to be known by context. \\

\begin{lemma} $\big[$\textit{Mean curvature comparison for $N$-Bakry-Émery Ricci tensor} $\big]$ \newline
\label{mean_curv_comp_N_bakry-lemma}
Let $(\textsl{M},\textsl{g}, e^{-f}d \textit{vol}_{\textsl{g}})$ be a SMMS and $\Sigma \subset \textsl{M}$ a hypersurface. If $\Riccfn(\dot{{\gamma}}, \dot{{\gamma}}) \geq 0$, along any minimal geodesic from $\Sigma$, then

\begin{flalign}
\label{mean_curv_comp_N_bakry}
\hspace{0.8cm}
 \emef(r) & \leq \emeh(r)  &  
\end{flalign}
holds up to the first occurrence of a cut-point.
\end{lemma}
\textit{Proof.}
Applying the corresponding Bochner inequality (\ref{bochner_f_N_bakry}) to the distance from $\Sigma$, we get 

\begin{flalign*}
\hspace{0.8cm}
0 = & \frac{1}{2} \Delta_{f}\vert \nabla r \vert^2 \geq \frac{\emef^{2}}{n+N-1} + \frac{\1}{\1 r} \emef. &
\end{flalign*}
Since $\emeh^{\prime}(r) = - \frac{1}{n+N-1}\emeh(r)^2$ we immediately get (\ref{mean_curv_comp_N_bakry}), by standard Sturm-Liouville comparison arguments, observing the initial conditions of $\emef$ and $\emeh$ match. $\hspace{0.3cm} \square$ \vspace{0.2cm}

\begin{lemma} $\big[$\textit{Mean curvature comparison for $\infty$-Bakry-Émery} $\big]$ \newline
\label{mean_curv_comp_infty_bakry-lemma}
Let $(\textsl{M},\textsl{g}, e^{-f}d \textit{vol}_{\textsl{g}})$ be a SMMS and $\Sigma \subset \textsl{M}$ a hypersurface. Suppose $\Riccf(\dot{{\gamma}}, \dot{{\gamma}}) \geq 0$, along any minimal geodesic from $\Sigma$. If $\frac{\1}{\1r}(f)\geq 0$ and $\hasf \geq 0$ along $\Sigma$ then

\begin{flalign}
\label{mean_curv_comp_infty_bakry}
\hspace{0.8cm}
 \emef(r) & \leq \emeh(r) &
\end{flalign}
holds up to the first occurrence of a cut-point.
\end{lemma}
\textit{Proof.}
By the usual Bochner formula applied to the distance (\ref{distance_to_hyper}), we have $0 \geq \eme^2/(n-1) + \eme^{\prime} + \Ricc(\nabla r, \nabla r)$. Summing and subtracting $D^{2}f(\nabla r, \nabla r)$ and using the curvature condition $\Riccf(\1_r, \1_r) \geq 0$, we obtain:

\begin{flalign*}
\hspace{0.8cm}
(\eme - \emeh)^{\prime}
	& \leq  -\frac{\eme^{2}}{n-1} -\Ricc (\nabla r, \nabla r) + \frac{\emeh^{2}}{n-1} \\
	& \leq -\frac{\eme^{2}- \emeh^{2}}{n-1} + D^{2}f(\nabla r, \nabla r). &
\end{flalign*}
Defining the auxiliary function $\hit(r) = (n-1) + \hasf(x) r$, $x \in \Sigma$ so that $(h^2)^{\prime} = \frac{2}{n-1}h^2 \emeh$ we henceforth obtain

\begin{flalign*}
\hspace{0.8cm}
\big[h^{2}(\eme - \emeh) \big]^{\prime}  
	& \leq h^{2}\bigg[ -\frac{(\eme - \emeh)^2}{n-1} +D^{2}f (\nabla r, \nabla r) \bigg] \\
	& \leq h^{2} D^{2}f (\nabla r,\nabla r). &
\end{flalign*}
Integrating this last inequality on $[0,r]$, and using integration by parts we get

\begin{flalign*}
\hspace{0.8cm}
h^{2}(r)\big(\eme(r) - \emeh(r)\big) 
	& - h^{2}(0)\big(\eme(0) - \emeh(0)\big)  \leq \int_{0}^{r}h^{2}(t)\frac{\1^{2} f}{\1 t^2}(t)dt \\
	& = h^{2}(r)\1_t f(r) - h^{2}(0)\1_t f(0) - \int_{0}^{r}\frac{\1 h(t)^{2}}{\1 t}\frac{\1 f}{\1 t}(t)dt &
\end{flalign*}
which is equivalent to

\begin{flalign*}
\hspace{0.8cm}
h^{2}(r)(\emef(r) - \emeh(r)) & \leq - \int_{0}^{r}(h(t)^{2})^{\prime}\1_t f(t)dt, &
\end{flalign*}
since $\big(\eme(0) - \emeh(0) - \1_t f(0)\big) = 0$. Now, if $ \frac{\1 f}{\1 r} \geq 0$ and if $\hasf (x)\geq 0$ then 

\begin{flalign*}
\hspace{0.8cm}
h^{2}(r)(\emef(r) - \emeh(r)) & \leq 0 &
\end{flalign*}
and the claim follows. $ \hspace{0.3cm} \square$ \vspace{0.2cm} 

It is now straightforward to get comparisons for the volume element based on lower bounds for the Bakry-Émery Ricci tensors (\ref{finite_bakry-emery_Ricci-T}) and (\ref{infinite_bakry-emery_Ricci-T}). Recall first that if $\mathcal{A}$ denotes the function on $[0,\infty) \times \Sigma$ defined by $(\Phi^{*}d \textit{vol}_{\textsl{g}})_{(r, x)} = \mathcal{A}(r,x) d\sigma \wedge dr$, for $x \in \Sigma$ and $0 \leq r < \tau(x)$, then $\frac{\1}{\1 t}\big\vert_{t=r} \lna \mathcal{A}(t,x) = \eme(r)$ so that if $\mathcal{A}_{f}(t,x) = e^{-(f \circ \Phi)(t,x)}\mathcal{A}(t,x)$ then

\begin{flalign*}
\hspace{0.8cm}
\frac{\1}{\1 t}\bigg\vert_{t=r} \lna \mathcal{A}_f(t,x) & = - df\Big( \frac{\1}{\1 r} \Big) + \eme(r) = \emef(r). &
\end{flalign*}
Define the function $\theta: [0,\infty) \times \Sigma \to \mathbb{R}$ by

\begin{flalign*}
\hspace{0.8cm}
\theta(t,x) &= \frac{\mathcal{A}_f(t, x)}{\Big( 1+ \frac{\hasf(x)}{k-1}t \Big)^{k-1}}, &
\end{flalign*}
where $k=n+N$ or $k = n$ depending on whether $N < \infty$ or $N = \infty$, respectively. It is immediate that, for each $x \in \Sigma$, we have

\begin{flalign*}
\hspace{0.8cm}
\frac{\1}{\1 t}\bigg\vert_{t=r} \theta(t,x) & = \theta (r,x) \Bigg( \emef(r) - \frac{\hasf(x)}{1 + \frac{\hasf(x)}{k-1} r} \Bigg). &
\end{flalign*}
Finally, on the one hand, if we have a nonnegative $N$-Bakry-Émery Ricci tensor then by Lemma \ref{mean_curv_comp_N_bakry-lemma}, we get $\frac{\1}{\1 t} \theta  = \theta \big\{ \emef - \emeh^{n+N} \big\} \leq 0$. On the other hand, if we have a nonnegative $\infty$-Bakry-Émery Ricci tensor and if $\frac{\1}{\1r}(f)\geq 0$ along minimal geodesics from $\Sigma$ and if $\hasf \geq 0$ on $\Sigma$ then by Lemma \ref{mean_curv_comp_infty_bakry-lemma} we have $\frac{\1}{\1 t} \theta = \theta \big\{ \emef - \emeh^n \big\} \leq 0$. In each case the function $\theta$ is monotone decreasing for each $x \in \Sigma$ so that $\theta(t,x) \leq \theta(0,x) = e^{-f(\Phi(0,x))}$. In any case

\begin{flalign}
\label{vol-elem_dominated}
\hspace{0.8cm}
\mathcal{A}_f(r, x) & \leq e^{-f(x)} \bigg( 1+ \frac{\hasf(x)}{k-1}r \bigg)^{k-1}, \hspace{0.2cm} x\in \Sigma, \hspace{0.2cm} 0\leq r< \tau(x), &
\end{flalign}
where $k=n+N$ or $k = n$ depending on whether we have $\Riccfn\geq 0$ or $\Riccf \geq 0$, respectively.

\subsection{The $f$-asymptotic volume ratio}
\label{f_asymptotic}

In this section we review in more detail the notion of asymptotic weighted volume ratio corresponding to the setting of a SMMS with nonnegative Bakry-Émery Ricci curvature, briefly described at the introduction. \vspace{0.2cm}

Suppose $(\textsl{M},\textsl{g})$ is a complete, noncompact, $n-$dimensional Riemannian manifold and let $\pzero \in \textsl{M}$ be an arbitrary (base) point. Let $\rho(w) = \text{dist}_{\textsl{g}}(w,\pzero)$, $w\in \textsl{M}$, be the distance to the point $\pzero$, and $f\in C^{\infty}(\textsl{M})$ be any function.

\begin{itemize}
\item If $\Riccfn \geq 0$ for some $N > 0$ then, it is well known (cf. \cite{qian1997estimates}, \cite{bakry2005volume}, \cite{wei2007comparison}) that the function 
\begin{flalign*}
\hspace{0.8cm}
(0, \infty) \ni R \mapsto & \Theta_{f}(R) \doteq \frac{\textit{vol}_{f}(B_{\pzero}(R))}{\w _{n+N} R^{n+N}} &
\end{flalign*}
is nonincreasing.
\item If $\Riccf \geq 0$ and $\frac{\1}{\1 \rho} (f) \geq 0$ along minimal geodesic segments from $\pzero$, then it is known (cf. \cite{wei2009comparison}) that the function  
\begin{flalign*}
\hspace{0.8cm}
(0, \infty) \ni R \mapsto & \Theta_f(R) \doteq \frac{\textit{vol}_{f}(B_{\pzero}(R))}{\w _{n} R^{n}} &
\end{flalign*}
is nonincreasing.
\end{itemize}
Now, regardless of the curvature assumption, there is a monotonous quantity denoted by the symbol $\Theta_{f}$. Henceforth, we convention that the quantity $\Theta_{f}$ we are referring to depend on which curvature assumption is in effect. Based on this, we define the \textit{asymptotic weighted volume ratio} of the triple $(\textsl{M}, \textsl{g},e^{-f}\textit{dvol}_{\textsl{g}} )$ as being its $f\textit{-}$asymptotic volume ratio, that is

\begin{flalign*}
\hspace{0.8cm}
f \textit{-}AVR(\textsl{g}) & = \lim_{r \rightarrow \infty} \Theta_{f}(r). &
\end{flalign*}

\begin{exmp} $\big[$\textit{Gaussian soliton} $\big]$ \newline
\label{gaussian_sol} Let $\textsl{M} \doteq \mathbb{R}^n$ with the usual Euclidean metric $\bar{g}$ and let $f(x) = (\lambda/2)\vert x \vert^2$, so that $\bar{D}^2f = \lambda \bar{g}$ and $\Riccf = \lambda \bar{g}$. This example shows that, although the underlying ambient is noncompact, we may have $\Riccf \geq \lambda \bar{g}$ with $\lambda >0$. Moreover

\begin{flalign*}
\hspace{0.8cm}
\lim_{r \rightarrow \infty } \textit{vol}_{f}(B_{\textbf{0}}(r)) & = \int_{\mathbb{R}^n} e^{-\frac{\lambda}{2}\vert x \vert^2}dx < \infty &
\end{flalign*}
and $ \langle \bar{\nabla} f, \bar{\nabla} \rho \rangle = \lambda \vert x \vert>0$, for $\rho(x)=\vert x\vert$, so that

\begin{flalign*}
\hspace{0.8cm}
f \textit{-}AVR(\bar{g}) & = 0. \hspace{0.3cm} \square \vspace{0.2cm}&
\end{flalign*}
\end{exmp}

\begin{theo} $\big[$\textit{Volume finiteness I} $\big]$ \newline
Let $(\textsl{M}, \textsl{g})$ be complete Riemannian manifold and $f \in C^{\infty}(\textsl{M})$ such that $\Riccf \geq \lambda >0$. Then $\textit{vol}_{f}(\textsl{M})$ is finite and $\textsl{M}$ has finite fundamental group. 
\end{theo}

For the $N < \infty$ case, we can reach stronger conclusions and draw comparisons with the classical Myers' Theorem (cf. \cite{wei2007comparison}). \vspace{0.2cm}

\begin{theo} $\big[$\textit{Volume finiteness II} $\big]$ \newline
Let $(\textsl{M}, \textsl{g})$ be a complete Riemannian manifold and let $f \in C^{\infty}(\textsl{M})$ be a function such that $\Riccfn \geq (n-1)H >0$, with $N < \infty$. Then, $\textsl{M}$ is compact and $\text{diam}_{\textsl{g}}(M) \leq \sqrt{\frac{n+N-1}{n-1}} \frac{\pi}{\sqrt{H}}.$
\end{theo}


\section{Proof of Theorems \ref{main_theo_N_bakry-emery} and \ref{main_theo_infty_bakry-emery}}
\label{proofs_section}
This section is divided into two parts. First, we give the proof of Theorem \ref{main_theo_N_bakry-emery} and, in the final subsection, we outline the proof of Theorem \ref{main_theo_infty_bakry-emery}. To prove both inequalities we may assume, without loss of generality, that $\Omega$ has no hole, that is: $\textsl{M}\setminus \Omega$ has no bounded component. As always, we denote by $\mathcal{T}_{\Omega}(R)$ the geodesic tube of radius $R$ about $\Omega$, i.e. the set $\mathcal{T}_{\Omega}(R) \doteq \{w \in \textsl{M} :\text{dist}_{\textsl{g}}(w, \Omega ) < R \}$.

\subsection{A sharp geometric inequality for $\Riccfn \geq 0$}
\label{prf_ineq_four}

\textbf{Proof of inequality (\ref{willmorex_type_N_bakry-emery}).} We compute the volume of $\mathcal{T}_{\Omega}(R)$ wrt the measure $e^{-f}d \textit{vol}_{\textsl{g}}$ to get

\begin{flalign*}
\hspace{0.8cm}\textit{vol}_{f}(\mathcal{T}_{\Omega}(R)) 
	& = \textit{vol}_{f}(\Omega) + \int_{\Sigma} \int_{0}^{R \wedge \tau(x)} e^{-f( \Phi(t,x))}\mathcal{A}(t,x) dt d\sigma (x) \notag \\ 
	& \leq  \textit{vol}_{f}(\Omega) +
\int_{\Sigma} \int_{0}^{R \wedge \tau(x)} e^{-f(x)}\bigg( 1+ \frac{\hasf(x)}{n+N-1}t \bigg)^{n+N-1} dt d\sigma (x) \notag \\ 
	& \leq  \textit{vol}_{f}(\Omega) +
\int_{\Sigma} \int_{0}^{R \wedge \tau(x)} e^{-f(x)}\bigg( 1+ \frac{\hasf^{+}(x)}{n+N-1}t \bigg)^{n+N-1} dt d\sigma (x) \notag \\ 
	& \leq \textit{vol}_{f}(\Omega) +
\int_{\Sigma} \int_{0}^{R} e^{-f(x)}\bigg( 1+ \frac{\hasf^{+}(x)}{n+N-1}t \bigg)^{n+N-1} dt d\sigma (x) \notag \\ 
	& = \textit{vol}_{f}(\Omega)+ \frac{R^{n+N}}{n+N} \int_{\Sigma} e^{-f(x)}\bigg(\frac{\hasf^{+}(x)}{n+N-1}\bigg)^{n+N-1} d\sigma(x) +\mathcal{O}(R^{n+N-1}),&
\end{flalign*}
where we have used inequality $(\ref{vol-elem_dominated})$ in the second line. Dividing both sides of the inequality above by $\omega_{n+N} R^{n+N} = \vert \mathbb{S}^{n+N-1} \vert R^{n+N}/(n+N)$ and letting $R \rightarrow \infty$ we obtain

\begin{flalign*}
\hspace{0.8cm}
f \textit{-}AVR(\textsl{g}) & \leq \frac{1}{\big\vert\mathbb{S}^{n+N-1} \big\vert} \int_{\Sigma} \bigg(\frac{\hasf^{+}(x)}{n+N-1}\bigg)^{n+N-1} e^{-f(x)} d\sigma(x), &
\end{flalign*}
which clearly implies (\ref{willmorex_type_N_bakry-emery}). \vspace{0.1cm}

\textbf{Equality discussion.} Suppose we have

\begin{flalign}
\label{equal_disc_N}
\hspace{0.8cm}
\int_{\Sigma} \bigg(\frac{\hasf^{+}(x)}{n+N-1}\bigg)^{n+N-1} e^{-f(x)} d\sigma(x) & = \big\vert\mathbb{S}^{n+N-1} \big\vert f \textit{-}AVR(\textsl{g}). &
\end{flalign}
First observe that from the continuity of $\tau$ and from the inequalities in \hyperref[prf_ineq_four]{Proof of inequality (5)} above we deduce $\tau \equiv \infty $ on the open set $\Sigma^{+} \doteq \{ x\in\Sigma : \hasf(x)>0 \}$. In addition, we have $\mathcal{A}(r, x) \leq e^{-f(x)}$ for all $x \in \Sigma \setminus \Sigma^+$, and all $0\leq r<\tau(x)$ by (\ref{vol-elem_dominated}). Therefore,

\begin{flalign*}
\hspace{0.8cm}\textit{vol}_{f}(\mathcal{T}_{\Omega}(R)) 
	& = \textit{vol}_{f}(\Omega) + \int_{\Sigma^{+}} \int_{0}^{R} \mathcal{A}_f(t,x) dt d\sigma (x)
	+ \int_{\Sigma \setminus \Sigma^{+}} \int_{0}^{R \wedge \tau_c(x)} \mathcal{A}_f(t,x) dt d\sigma (x) \notag \\
	& \leq \textit{vol}_{f}(\Omega) + \int_{\Sigma^{+}} \int_{0}^{R} \theta(t, x)  \bigg( 1+ \frac{\hasf(x)t}{n+N-1} \bigg)^{n+N-1} dt d\sigma (x)	+ \int_{\Sigma \setminus \Sigma^{+}} \int_{0}^{R} e^{-f(x)} dt d\sigma (x) \notag \\	
	& \leq \textit{vol}_{f}(\Omega) + \int_{\Sigma^{+}} \int_{R^{\prime}}^{R} \theta(t, x)  \bigg( 1+ \frac{\hasf(x)t}{n+N-1} \bigg)^{n+N-1} dt d\sigma (x)  \notag \\
	&  + \int_{\Sigma^{+}} \int_{0}^{R^{\prime}} \theta(t, x)  \bigg( 1+ \frac{\hasf(x)t}{n+N-1} \bigg)^{n+N-1} dt d\sigma (x)	+ \mathcal{O}(R) \\
	& \leq \textit{vol}_{f}(\Omega) + \int_{\Sigma^{+}} \theta(R^{\prime}, x)\int_{R^{\prime}}^{R}   \bigg( 1+ \frac{\hasf(x)t}{n+N-1} \bigg)^{n+N-1} dt d\sigma (x) \notag \\
	& + \int_{\Sigma^{+}} \int_{0}^{R^{\prime}} \theta(t, x)  \bigg( 1+ \frac{\hasf(x)t}{n+N-1} \bigg)^{n+N-1} dt d\sigma (x)	+ \mathcal{O}(R)	
\end{flalign*}
Dividing both sides by $\omega_{n+N} R^{n+N} = \vert \mathbb{S}^{n+N-1} \vert R^{n+N}/(n+N)$ and letting $R\rightarrow \infty$, we get

\begin{flalign*}
\hspace{0.8cm}
f \textit{-}AVR(\textsl{g}) & \leq \frac{1}{\big\vert\mathbb{S}^{n+N-1} \big\vert}
\int_{\Sigma^{+}} \bigg(\frac{\hasf^{+}(x)}{n+N-1}\bigg)^{n+N-1}  \theta(R^{\prime}, x) d\sigma(x). &
\end{flalign*}
Letting $R^{\prime} \rightarrow \infty$ yields

\begin{flalign*}
\hspace{0.8cm}
f \textit{-}AVR(\textsl{g}) & \leq \frac{1}{\big\vert\mathbb{S}^{n+N-1} \big\vert}
\int_{\Sigma^{+}} \bigg(\frac{\hasf^{+}(x)}{n+N-1}\bigg)^{n+N-1} \theta_{\infty}(x) d\sigma(x), &
\end{flalign*}
where $\theta_{\infty}(x) = \lim_{r \rightarrow \infty }\theta(t,x)$. As we have equality in (\ref{equal_disc_N}), and since $t\mapsto \theta (t,x)$ is nonincreasing, we ought to have $\theta_{\infty}(x) = e^{-f(x)}$ for a. e. $x \in \Sigma^{+}$, which implies

\begin{flalign*}
\hspace{0.8cm}
\mathcal{A}_f(t, x) & = e^{-f(x)}\Big( 1+ \frac{\hasf(x)}{n+N-1}t \Big)^{n+N-1} &
\end{flalign*}
for all $x\in \Sigma^+$ and all $t \geq 0$, by continuity. Differentiating $\mathcal{A}_f(t, x)$ with respect to $t$ immediately gives

\begin{flalign}
\label{equality_N-laplacian}
\hspace{0.8cm}
\emef(r) & = \frac{\1}{\1 t}\bigg\vert_{t=r} \lna \mathcal{A}_f(t, x)  = \frac{(n+N-1)\hasf(x) }{n+N-1 + \hasf(x) r} = \emeh (r), \hspace{0.2cm} x \in \Sigma^+. &
\end{flalign}
Inspecting the corresponding Bochner formula (\ref{bochner_f_N_bakry}) on $\Phi ([0, \infty)\times \Sigma^+)$, we must have

\begin{enumerate}[label=\textit{(E\hspace{-0.035cm}\arabic*)},leftmargin=1.7cm]
\item $D^2 r = \frac{\Delta r}{n-1} g$;

\item $\Ricc_f^N(\1_r,\1_r)=0$;

\item $\Delta r = -\frac{n-1}{N}\frac{\1}{\1 r}(f)$.
\end{enumerate}
Since $\Ricc_f^N \geq 0$, it follows from $\textit{(E\hspace{-0.035cm}2)}$ that $\Ricc_f^N (\1_r, \cdot ) = 0$. It follows from $\textit{(E\hspace{-0.035cm}1)}$ that $\Sigma^+$ is umbilic, that is, $A_{\xi} = \frac{H}{n-1}Id$ where $Id$ is the identity on $T\Sigma^+$, and $\xi = -\grad r$. Employing an orthonormal frame $\{e_0=\nu, e_1,..., e_{n-1}\}$ along $\Sigma^+$, it follows from $\textit{(E\hspace{-0.035cm}3)}$ that 

\begin{flalign*}
\hspace{0.8cm}
D^2 f(e_j, \nu)& = e_j \langle \nabla f, \nu \rangle - \langle \nabla f, \nabla_{e_j} \1 r \rangle = - \frac{N \has_j}{n-1}-\frac{\has f_j}{n-1}. &
\end{flalign*}
Using the Codazzi equation, with $1 \leq j,k \leq n-1$, we get $\big(\textsl{R}(e_k, e_j)\xi\big)^{T} = \frac{\has_j}{(n-1)} e_k - \frac{\has_k}{(n-1)} e_j$ so that the Ricci curvature of the pair $(e_j,-\nu)$ is $\Ricc(e_j,\xi) = \frac{n-2}{n-1}\has_j$. Now, by definition (\ref{finite_bakry-emery_Ricci-T}) and equation $\textit{(E\hspace{-0.035cm}3)}$ again, we have

\begin{flalign*}
\hspace{0.8cm}
0 & = \Ricc_f^N(\nu, e_j) = -\frac{n-2}{n-1}\has_j - \bigg(\frac{N\has_j}{n-1}+\frac{\has f_j}{n-1}\bigg) - \frac{1}{N}\bigg( -\frac{N\has f_j}{n-1} \bigg) = -\frac{n+N-2}{n-1}H_j &
\end{flalign*}
whence $\has$ is locally constant on $\Sigma^+$. Moreover

\begin{flalign*}
\hspace{0.8cm}
\has & = -\frac{n-1}{N} \nu(f), \hspace{0.2cm} \text{and} \hspace{0.2cm} \hasf  = \has - \bigg( -\frac{N\has}{n-1} \bigg) = \frac{n+N-1}{n-1}\has,   &
\end{flalign*}
holds on $\Sigma^+$ by equation $\textit{(E\hspace{-0.035cm}3)}$ and, as a result, both $\hasf$ and the normal derivative $\nu(f)$ are also locally constant where $\has$ is. We conclude that $\Sigma^+$ is the union of several connected components of $\Sigma$. Since $\Phi$ is a diffeomorphism from $[0, \infty)\times \Sigma^+$ onto its image, we have $\Phi^*g = dr\otimes dr + \hat{g}_r$ where $\hat{g}_r$ is a family of metrics on $\Sigma^+$ with $\hat{g}_0 = g_{\Sigma^+}$. In what follows, we integrate the radial derivatives $\frac{\1}{\1 r}(f)$ along minimal geodesics $[0, \infty) \ni t\mapsto \exp_x t\nu_x$, from $x \in \Sigma^+$. Given local coordinates $(x^1,...,x^{n-1})$ on $\Sigma^+$ it follows from $\textit{(E\hspace{-0.035cm}3)}$ that 

\begin{flalign*}
\hspace{0.8cm}
\frac{1}{2}\frac{\1}{\1 r} (\hat{g}_r)_{ij} & = -\frac{1}{N}\frac{\1 f}{\1 r}(\hat{g}_r)_{ij}&
\end{flalign*}
which implies $\hat{g}_r = c(x)e^{-\frac{2}{N} f[\exp_x r\nu_x]}g_{\Sigma^+}$, where $c(x) = e^{\frac{2}{N}f(x)}$. In addition to that, by the equations (\ref{equality_N-laplacian}) and $\textit{(E\hspace{-0.035cm}3)}$ and because $\hasf  = \frac{n+N-1}{n-1}\has$ on $\Sigma^+$, we have

\begin{flalign*}
\hspace{0.8cm}
\frac{\1 f}{\1 r} & = - \frac{N\has}{n-1 + \has r}. &
\end{flalign*}
Integrating this identity on $[0,r]$, gives

\begin{flalign*}
\hspace{0.8cm}
f\big[\exp_x r\nu(x)\big] & = \lna\bigg( 1+ \frac{\has r}{n-1}\bigg)^{-N} + f(x), \hspace{0.2cm}\text{for} \hspace{0.2cm} x \in \Sigma^+. &
\end{flalign*}
This results in writing the nontrivial part of the pulled-back metric on $[0, \infty)\times \Sigma^+$ as $\hat{g}_r = \big( 1+ \frac{\has r}{n-1}\big)^2 g_{\Sigma^+}$ and it also proves that $\Phi ([0, \infty)\times \Sigma^+)$ is isometric to $\big([r_0, \infty )\times \Sigma^+, dr^2 + ( r/r_0)^2 g_{\Sigma^+}\big)$, where $r_0 = \frac{n-1}{\has}$. If $\Sigma$ is connected and if we have the stronger equality

\begin{flalign*}
\hspace{0.8cm}
\int_{\Sigma}\bigg \vert \frac{ \hasf(x)}{n+N-1} \bigg\vert^{n+N-1} e^{-f(x)}d\sigma (x) & = \big\vert\mathbb{S}^{n+N-1} \big\vert f \textit{-}AVR(\textsl{g}), &
\end{flalign*}
then the previous analysis goes through to obtain that $\Sigma$ has positive constant weighted mean curvature and therefore $\Sigma^+=\Sigma$. The rest of the proof can be taken from \cite{wu2025willmore}. $\hspace{0.3cm} \square$

\subsection{A sharp geometric inequality for $\Riccf \geq 0$}
\label{prf_ineq_eight}

\textbf{Proof of inequality (\ref{willmorex_type_infty_bakry-emery}).} We have

\begin{flalign*}
\hspace{0.8cm}\textit{vol}_{f}(\mathcal{T}_{\Omega}(R)) 
	& = \textit{vol}_{f}(\Omega) + \int_{\Sigma} \int_{0}^{R \wedge \tau(x)} \mathcal{A}_f(t,x) dt d\sigma (x) \notag \\ 
	& \leq  \textit{vol}_{f}(\Omega) + \int_{\Sigma} \int_{0}^{R \wedge \tau(x)} e^{-f(x)}\bigg( 1+ \frac{\hasf(x)}{n-1}t \bigg)^{n-1} dt d\sigma (x) \notag \\ 
	& \leq	\textit{vol}_{f}(\Omega) +
\int_{\Sigma} \int_{0}^{R} e^{-f(x)}\bigg( 1+ \frac{\hasf(x)}{n-1}t \bigg)^{n-1} dt d\sigma (x) \notag \\
	& = \textit{vol}_{f}(\Omega)+ \frac{R^{n}}{n} \int_{\Sigma} e^{-f(x)}\bigg(\frac{\hasf(x)}{n-1}\bigg)^{n-1} d\sigma(x) +\mathcal{O}(R^{n-1}),&
\end{flalign*}
Dividing both sides of the inequality above by $\omega_{n} R^{n} = \vert \mathbb{S}^{n-1} \vert R^{n}/n$ and letting $R \rightarrow \infty$ we obtain (\ref{willmorex_type_infty_bakry-emery}).

\textbf{Equality discussion.} Suppose we have

\begin{flalign*}
\hspace{0.8cm}
\int_{\Sigma}\bigg( \frac{ \hasf(x)}{n-1} \bigg)^{n-1} e^{-f(x)}d\sigma (x) & = \vert\mathbb{S}^{n-1} \vert f \textit{-}AVR(\textsl{g})  &
\end{flalign*}
Again, $\tau \equiv \infty $ on $\Sigma^{+} \doteq \{\hasf(x)>0 \}$ and by exactly the same argument as before, $\theta_{\infty}(x) = e^{-f(x)}$ for a. e. $x \in \Sigma^+$ which implies $\mathcal{A}_f(t, x) = e^{-f(x)}\big( 1+\frac{H(x)}{n-1}t \big)^{n-1}$, for all $x \in \Sigma^+$ and all $t\geq 0$. Differentiating  $t \mapsto \mathcal{A}_f(t,x)$ yields

\begin{flalign*}
\hspace{0.8cm}
\emef(r) & = \frac{\1}{\1 t}\bigg\vert_{t=r}\lna \mathcal{A}_f(t, x) = \emeh(r).&
\end{flalign*}
for all $r \geq 0$, and $x\in \Sigma^+$. From Lemma \ref{mean_curv_comp_infty_bakry-lemma}, we deduce  $\1 f/\1 r \equiv 0$ so that, outside $\Omega$, $f$ is constant along all minimal geodesic segments from $\Sigma^+$. Looking into the Bochner inequality from (\ref{mean_curv_infty-bakry}) on $\Phi ([0, \infty)\times \Sigma^+)$, we get

\begin{enumerate}[label=\textit{(E\hspace{-0.035cm}\arabic*)},leftmargin=1.7cm]
\item $D^2 r = \frac{\Delta r}{n-1} g$;

\item $\Ricc_f(\1_r,\1_r)=0$.
\end{enumerate}
It follows from $\textit{(E\hspace{-0.035cm}1)}$ that $\Sigma^+$ is umbilic. Again, taking an orthonormal frame $\{ \nu=e_0, e_1,..., e_{n-1}\}$ on $\Sigma^+$, it follows from $\textit{(E\hspace{-0.035cm}2)}$ together with $\Riccf\geq 0$ that $\Riccf(\1 r, \cdot )=0$. As in the previous rigidity argument, it follows from the Codazzi equation that

\begin{flalign*}
\hspace{0.8cm}
0 & = \Ricc_f(\nu, e_j) = \big(\Ricc + D^2f \big)(\nu, e_j) = \frac{n-2}{n-1}\has_j + \frac{\has}{n-1}f_j, & 
\end{flalign*}
thus $f = \lna \has^{n-2}$ on $\Sigma^+$. Now, because $\Phi:[0, \infty)\times \Sigma^+\to \Phi([0, \infty)\times \Sigma^+)$ is a diffeomorphism, the pulled-back metric takes the form $dr^2+\hat{g}_r$, for a certain one-parameter family of metrics $\hat{g}_r$ on $\Sigma^+$ such that $\hat{g}_0 = g_{\Sigma^+}$. If $(x^1,...,x^{n-1})$ are local coordinates on $\Sigma^+$, then by $\textit{(E\hspace{-0.035cm}1)}$ we have

\begin{flalign*}
\hspace{0.8cm}
\frac{1}{2}\frac{\1}{\1 r} (\hat{g}_r)_{jk} & = \frac{\has}{n-1 + \has r } (\hat{g}_r)_{jk}, \hspace{0.2cm}\text{for} \hspace{0.2cm} 1\leq j,k\leq n-1 &
\end{flalign*} 
whence $\hat{g}_r = \big( 1+ \frac{\has(x)}{n-1}r \big)^2g_{\Sigma^+}$. This means that the portion of $\textsl{M}\setminus \Omega$ which is the image of $\textsl{N}\Sigma^+$  under the normal exponential map is isometric to the twisted product metric

\begin{flalign*}
\hspace{0.8cm}
\bigg( [0, \infty)\times \Sigma^+,& dr^2 + \bigg( 1+\frac{\has(x)}{n-1}r  \bigg)^2g_{\Sigma^+}\bigg).&
\end{flalign*}
If $\has >0$ over $\Sigma$ then $\Sigma^+ = \Sigma$ and the previous argument applies to the whole of $\Sigma$. In particular, by compactness of $\Sigma$, the mean curvature $\has$ is bounded above and below by positive constants. Therefore $f = \lna \has^{n-2}$ is bounded. Moreover, if $\textsl{M}$ has only one end then $\Sigma$ is connected. $\hspace{0.3cm} \square$

\section{Acknowledgements}
The author is grateful to Professor Ezequiel Barbosa for suggesting the topic during the author's PhD.

%
%
%

\bibliographystyle{acm}
\bibliography{references}
\end{document}